\documentclass[12pt]{article}
\usepackage{amssymb,amsmath,amsthm,tabularx}

\theoremstyle{plain}
\newtheorem{thm}{Theorem}[section]
\newtheorem{cor}[thm]{Corollary}

\theoremstyle{definition}
\newtheorem{rem}[thm]{Remark}

\newcommand{\bb}[1]{\mathbb{#1}}
\newcommand{\cl}[1]{\mathcal{#1}}
\newcommand{\n}{\noindent}

\newcommand{\vp}{\varepsilon}

\begin{document}

\title{APPROXIMATION PROPERTIES FOR CROSSED PRODUCTS BY ACTIONS AND COACTIONS}
\author{May M. Nilsen \hspace{.7in} and \hspace{.7in} Roger R.
Smith\footnote{Partially supported by  a grant from the National Science
Foundation.}\\ {\tt mnilsen@math.tamu.edu\hspace{1.2in}
rsmith@math.tamu.edu}\\ {}\\{}\\Department of Mathematics\\
Texas A\&M University\\
College Station, TX \ 77843-3368\\{}\\{}\\{}}

\date{}
\maketitle

\begin{abstract}We investigate approximation properties for $C^*$-algebras and
their crossed products by actions and coactions by locally compact groups. We
show that Haagerup's approximation constant is preserved for crossed products
by arbitrary amenable groups, and we show why this is not always true in
the non-amenable case. We also examine similar questions for other forms of the
approximation property. \end{abstract}\vspace{1in}
\n 2000 Mathematics Subject Classification Numbers: 46L55, 46L07, 46L06
\thispagestyle{empty}
\newpage
\pagestyle{plain}
\setcounter{page}{1}

\baselineskip = 19.05pt

\section{Introduction}\label{sec1}

\indent

A $C^*$-algebra $\cl{A}$ has the {\it  completely bounded approximation
property} (CBAP) if there is a net $\{T_\gamma\colon \ \cl{A}\to
\cl{A}\}_{\gamma\in\Gamma}$ of finite rank maps, uniformly bounded in the
completely bounded norm, which converges in the point norm topology to the
identity. The smallest number which can bound such a net is called the
Haagerup constant of $\cl{A}$, and is denoted by $\Lambda(\cl{A})$. If no such
net exists we set $\Lambda(\cl{A})$ equal to $\infty$. This constant was
introduced and studied in a series of papers \cite{CH,DH,H}, and is an
important isomorphism invariant for $C^*$-algebras. An interesting problem is
to determine the functorial properties of the CBAP. In \cite{SS1}, it was
shown that $\Lambda$ is multiplicative on the minimal tensor product
$\cl{A}\otimes \cl{B}$ of $C^*$-algebras, and $\Lambda$ is invariant under
crossed products by discrete amenable groups, \cite{SS2}. The main result of
the paper is the extension to the case of general locally compact amenable
groups. For von~Neumann algebras $\cl{M}$ there is a corresponding constant
$\Lambda_w(\cl{M})$, where point norm convergence is replaced by point
$w^*$-convergence. The equality \begin{equation}\label{eq1.1}
\Lambda_w(\cl{M}\times_\alpha G) = \Lambda_w(\cl{M}) \end{equation}
for amenable $W^*$-dynamical systems $(\cl{M},G,\alpha)$ was obtained in \cite{HK}
and independently in \cite{A}. The difference between norm and
$w^*$-convergence means that techniques appropriate for $\Lambda_w$ rarely
carry over to $\Lambda$, and this is the case here. It also seems impossible to
extend the methods of \cite{SS2} to general amenable groups.

Our approach to this problem will make heavy use of the duality theorems of
Katayama, \cite{Ka}, and of Imai and Takai, \cite{IT}. This requires us to
study coactions, and so we have included in the second section a brief review
of those parts of the theory which we will use subsequently. The method which
underlies our work is to obtain information about a $C^*$-algebra $\cl{A}$ by
constructing approximate point norm factorizations of $\cl{A}$ through a second
$C^*$-algebra $\cl{B}$. There are many instances of this in the literature (see,
for example, \cite{CE,Ki}) so we have formalized this by calling it the
{\it completely contractive factorization property} (CCFP), defined in the
third section. When a pair $(\cl{A},\cl{B})$ has the CCFP, we show in
Theorem~\ref{thm3.1} that slice map properties (see \cite{Kr1,Kr2,W}) and
various approximation properties, including the CBAP, pass from $\cl{B}$ to
$\cl{A}$. We prove that the pair $(\cl{A},\cl{A}\times_{\alpha,r}G)$ has the
CCFP (Theorem~\ref{thm3.2}), from which the inequality
\begin{equation}\label{eq1.2} \Lambda(\cl{A}) \le
\Lambda(\cl{A}\times_{\alpha,r}G) \end{equation}
follows immediately for any locally compact group $G$. The reverse inequality,
for $G$ amenable, is obtained in the fourth section, by showing that the pair
$(\cl{A}\times_{\alpha,r}G, (\cl{A}\times_{\alpha,r}G) \times_{\hat\alpha}G)$ has the
CCFP for the dual coaction $\hat\alpha$. Our methods also give preservation of
other approximation properties by amenable crossed products by actions and
coactions (Theorem~\ref{thm4.6}). The last section contains some concluding
remarks. We indicate why our results cannot be extended beyond the amenable
case, and we take the opportunity to point out that some recent work of Ozawa,
\cite{O}, can be used to settle two open problems from the paper of Haagerup
and Kraus, \cite{HK}.

\section{Preliminaries}\label{sec2}

\setcounter{equation}{0}

\indent

A triple $(\cl{A},G,\alpha)$, where $\cl{A}$ is a $C^*$-algebra, $G$ is a locally
compact group, and $\alpha\colon \ G\to \text{Aut}(\cl{A})$ is a homomorphism, is
said to be a $C^*$-dynamical system if the map $t\to \alpha_t(a)$ is
norm-continuous on $G$ for each $a\in \cl{A}$. The reduced crossed product
$\cl{A}\times_{\alpha, r}G$ is constructed by taking a faithful representation
$\pi\colon \ \cl{A}\to B(H)$ and associating to it a representation
$\tilde\pi\colon \ \cl{A}\to B(L^2(G,H))$ defined by
\begin{equation}\label{eq2.1}
(\tilde\pi(a)\xi)(t) = \pi(\alpha_{t^{-1}}(a))\xi(t)
\end{equation}
for $\xi\in L^2(G,H)$. Each $s\in G$ corresponds to a unitary $\lambda_s \in
B(L^2(G,H))$ defined by
\begin{equation}\label{eq2.2}
(\lambda_s\xi)(t) = \xi(s^{-1}t), \qquad \xi\in L^2(G,H),
\end{equation}
and it is easily checked that
\begin{equation}\label{eq2.3}
\lambda_s\pi(a)\lambda^*_s = \pi(\alpha_s(a))
\end{equation}
for $s\in G, a\in \cl{A}$ (the pair $(\pi,\lambda)$ is called covariant). The
reduced crossed product is then the norm closure of the set of operators of
the form
\begin{equation}\label{eq2.4}
\int f(s) \tilde\pi(a) \lambda_s \ ds
\end{equation}
where $a\in \cl{A}, f\in K(G)$, the algebra of continuous functions of compact
support on $G$, and $ds$ is a fixed choice of left Haar measure on $G$. This
$C^*$-algebra is independent of the choice of $\pi$, since any covariant pair
$(\pi,\lambda)$, where $\pi$ is faithful, induces a faithful representation
$\tilde\pi \times \lambda$ of $\cl{A}\times_{\alpha, r}G$ (see Chapter~7 of
\cite{P}).

We now briefly review the definition of a reduced coaction. A non-degenerate
injective $*$-homomorphism $\delta$ from $\cl{A}$ into the multiplier algebra
$M(\cl{A}\otimes C^*_r(G))$ of $\cl{A}\otimes C^*_r(G)$ (where $\otimes$ always denotes
the minimal $C^*$-tensor product) is called a coaction if it satisfies
\begin{itemize}\label{eq2.5}
\item[(1)]\hspace{2in} $\displaystyle (\delta\otimes I)\delta = (I \otimes
\delta_G)\delta$\hfill (2.5)
\end{itemize}
where $\delta_G\colon \ C^*_r(G) \to M(C^*_r(G) \otimes C^*_r(G))$ is the
integrated form of the map $s\to s\otimes s$, $s\in G$; 
\begin{itemize}\label{eq2.6}
\item[(2)]\hspace{2in} $\displaystyle\delta(\cl{A})(I\otimes C^*_r(G)) \subseteq \cl{A}\otimes
C^*_r(G)$,\hfill (2.6)
\end{itemize}
(see Definition~1.1 of \cite{GL}). The first condition is called the coaction
identity. In the earlier literature, the second condition was called
non-degeneracy of the coaction.
\medskip

\addtocounter{equation}{2}

\n A triple $(\cl{A},G,\delta)$, where $\delta$ is a coaction, is called a
cosystem. The crossed product $\cl{A}\times_\delta G$ is constructed on $H\otimes
L^2(G)$ in a manner similar to crossed products by actions. For actions the
crossed product is generated by a copy of $\cl{A}$ and a copy of $C^*_r(G)$; for
coactions $C^*_r(G)$ is replaced by $C_0(G)$. The crossed product
$\cl{A}\times_\delta G$ is the norm closed span of the set of elements
$$\{(\pi\otimes I)(\delta(a)) (I\otimes M_f)\colon \ a\in \cl{A}, f\in C_0(G)\},$$
where $\pi$ is a faithful representation of $\cl{A}$ on $H$, and $M_f$ denotes
multiplication by $f$ on $L^2(G)$ (Definition~1.4 of \cite{GL}). The coaction
identity (2.5) and non-degeneracy (2.6) ensure  that this is a $C^*$-algebra,
\cite{R}. We also note that the crossed product is independent of the choice
of faithful representation. This follows from the work of Raeburn, \cite{R},
who showed that the full crossed product, defined by universal properties, is
isomorphic to the reduced crossed product that we have defined here. The papers
\cite{GL,L,NT,N,R} are good references for background material.

Various forms of approximation properties will appear subsequently, so we
review them here for the reader's convenience. We say that $(\cl{A},\cl{B},\cl{E})$ has the
{\it slice map property}, where $\cl{E}$ is a closed subspace of a
$C^*$-algebra $\cl{B}$, if any element $x\in \cl{A}\otimes \cl{B}$, whose
right slices by $\phi\in \cl{A}^*$ lie in $\cl{E}$, must be an element of
$\cl{A}\otimes \cl{E}$. Then $\cl{A}$ has the {\it general slice map property}
if every such triple $(\cl{A},\cl{B},\cl{E})$ has the slice map property. If we
restrict $\cl{B}$ to being the algebra of compact operators $K(H)$ on a
separable Hilbert space, then we say that $\cl{A}$ has the {\it slice map
property for subspaces of} ${K(H)}$. If there exists a net of finite
rank operators $\{T_\gamma\colon \ \cl{A}\to \cl{A}\}_{\gamma\in \Gamma}$ such
that $T_\gamma\otimes I$ converges to $I\otimes I$ in the point norm topology
on $\cl{A}\otimes K(\ell^2)$, then $\cl{A}$ is said to have the operator
approximation property (OAP), \cite{ER,Kr2}. If the same conclusion holds true
when $K(\ell^2)$ is replaced by any $C^*$-algebra, then we say that $\cl{A}$
has the {\it strong operator approximation property} (strong OAP). In general,
these two forms of the OAP are distinct, although they coincide for the reduced
$C^*$-algebras of discrete groups, \cite{HK}.

\section{The completely contractive factorization property}\label{sec3}

\setcounter{equation}{0}

\indent

We say that an ordered pair $(\cl{A},\cl{B})$ of $C^*$-algebras has the
{\it completely contractive factorization property} (CCFP) if the following
condition is satisfied. Given $\vp>0$, and $a_1,\ldots, a_n\in \cl{A}$, there
exist completely contractive maps $S\colon \ \cl{A}\to \cl{B}$ and $T\colon \
\cl{B}\to \cl{A}$ such that \begin{equation}\label{eq3.1}
\|T(S(a_i)) - a_i\| < \vp,\qquad 1\le i \le n.
\end{equation}
Since the set of pairs $(F,\vp)$ of finite subsets $F$ of $\cl{A}$ and positive
numbers $\vp$ can be given a partial order by $(F_1,\vp_1) \le (F_2,\vp_2)$ if
and only if $F_1\subseteq F_2$ and $\vp_2\le \vp_1$, it is clear that our
definition is equivalent to the existence of nets $\{S_\gamma\colon \ \cl{A}\to \cl{B},
\quad T_\gamma\colon \ \cl{B}\to \cl{A}\}_{\gamma\in\Gamma}$ of complete contractions
satisfying
\begin{equation}\label{eq3.2}
\lim_\gamma\|T_\gamma(S_\gamma(a))-a\|=0
\end{equation}
for all $a\in \cl{A}$. The point of introducing this concept is to have a simple
method of transferring properties from $\cl{B}$ to $\cl{A}$. This is formalized in the
following result, which will allow us subsequently to concentrate on showing
that a particular pair of $C^*$-algebras has the CCFP.

\begin{thm}\label{thm3.1}
Let $(\cl{A},\cl{B})$ be a pair of $C^*$-algebras having the CCFP.
\begin{itemize}
\item[(1)] If $\cl{C}$ is any $C^*$-algebra then the pairs $(\cl{A},\cl{B}\otimes \cl{C})$ and
$(\cl{A}\otimes \cl{C}, \cl{B}\otimes \cl{C})$ have the CCFP.
\item[(2)] If $(\cl{B},\cl{C})$ has the CCFP then so too does $(\cl{A},\cl{C})$.
\item[(3)] The $C^*$-algebra $\cl{A}$ inherits each of the following properties
from $\cl{B}$:
\begin{itemize}
\item[(i)] The completely bounded approximation property, and $\Lambda(\cl{A}) \le
\Lambda(\cl{B})$;
\item[(ii)] Nuclearity;
\item[(iii)] The general slice map property;
\item[(iv)] The slice map property for subspaces of $K(H)$;
\item[(v)] The operator approximation property;
\item[(vi)] The strong operator approximation property;
\item[(vii)] Exactness.
\end{itemize}
\end{itemize}
\end{thm}

\begin{proof}
(1)~~Fix an element $c_0\in \cl{C}$ of unit norm, and pick a $\phi\in
{\cl{C}}^*$ such that $\|\phi\| = \phi(c_0)=1$. If $\{S_\gamma\colon \
\cl{A}\to \cl{B},\quad T_\gamma\colon \ \cl{B}\to \cl{A}\}_{\gamma\in \Gamma}$
are nets of complete contractions such that $\lim\limits_\gamma T_\gamma
S_\gamma = I$ in the point norm topology, define $S'_\gamma\colon \ \cl{A}\to
\cl{B}\otimes \cl{C}$, $T'_\gamma\colon \ \cl{B}\otimes \cl{C} \to \cl{A}$ by
\begin{equation}\label{3.2a}S'_\gamma(a) = S_\gamma(a)\otimes c_0,\quad
T'_\gamma(b\otimes c) = \phi(c) T_\gamma(b),\end{equation}
for $a\in \cl{A}, b\in \cl{B}, c\in \cl{C}, \gamma\in \Gamma$. These maps are complete
contractions, and $\lim\limits_\gamma T'_\gamma S'_\gamma = I$ in the point
norm topology, showing that $(\cl{A}, \cl{B}\otimes \cl{C})$ has the CCFP. To show that
$(\cl{A}\otimes \cl{C}, \cl{B}\otimes \cl{C})$ has the CCFP, just define $S''_\gamma = S_\gamma
\otimes I$ and $T''_\gamma = T_\gamma \otimes I$.

\n (2)~~This is a simple exercise in the composition of maps.

\n (3) (i)~~If $\cl{B}$ has the CBAP then there exists a net $\{R_\mu\colon \ \cl{B}\to
\cl{B}\}_{\mu\in M}$ of finite rank maps converging in the point norm topology to
$I$ and satisfying $\|R_\mu\|_{cb}\le \Lambda(\cl{B})$. Given $a_1,\ldots, a_n \in
\cl{A}$  and $\vp>0$, we can select $\gamma\in\Gamma$ and $\mu\in M$ such that
\begin{equation}\label{eq3.3}
\|T_\gamma R_\mu S_\gamma(a_i) -a_i\| < \vp, \qquad 1\le i \le n,
\end{equation}
and the result follows, since $\|T_\gamma R_\mu S_\gamma\|_{cb} \le \Lambda(\cl{B})$
and all such maps are finite rank.

\n (ii)~~It was shown in \cite{S} that nuclearity of $\cl{B}$ is equivalent to
having an approximate point norm factorization of $I$ through matrix algebras
by complete contractions. This in turn is clearly equivalent to
$(\cl{B},K(\ell^2))$ having the CCFP, so the result follows from the transitivity
of (2), taking $\cl{C}$ to be $K(\ell^2)$.

\n (iii)~~Suppose that $\cl{B}$ has the general slice map property. Let $\cl{C}$ be a
$C^*$-algebra with a closed subspace $\cl{E}$ and consider an element  $x\in
\cl{A}\otimes \cl{C}$, all of whose right slices lie in $\cl{E}$. If $\phi\in \cl{B}^*$, then the
composition of $S_\gamma\otimes I\colon \ \cl{A}\otimes \cl{C} \to \cl{B}\otimes \cl{C}$ with
$R_\phi$ is a right slice map on $\cl{A}\otimes \cl{C}$. By hypothesis $R_\phi((S_\gamma
\otimes I)(x))\in \cl{E}$, so all right slices of $(S_\gamma\otimes I)(x)$ lie in
$\cl{E}$. Since $\cl{B}$ has the general slice map property, $S_\gamma \otimes I(x) \in
\cl{B}\otimes \cl{E}$, so $T_\gamma S_\gamma \otimes I(x) \in \cl{A} \otimes \cl{E}$ for all
$\gamma\in \Gamma$. Take a limit over $\Gamma$ to show that $x\in \cl{A}\otimes \cl{E}$,
proving that $\cl{A}$ has the general slice map property.

\n (iv)~~This is a special case of (iii).

\n (v),(vi)~~We prove only (vi) since the argument also applies to (v).
Let $\cl{C}$ be any $C^*$-algebra, and suppose that $\cl{B}$ has the strong OAP. Given
$x_1,\ldots, x_n \in \cl{A}\otimes \cl{C}$, we may choose $\gamma\in\Gamma$ so that
\begin{equation}\label{eq3.4}
\|T_\gamma S_\gamma \otimes I(x_i) - x_i\| < \vp/2, \qquad 1\le i \le n.
\end{equation}
Then we may choose, by hypothesis, a finite rank map $R\colon \ \cl{B}\to \cl{B}$ such
that
\begin{equation}\label{eq3.5}
\|R\otimes I(S_\gamma\otimes I)(x_i) - S_\gamma\otimes I(x_i)\| < \vp/2,
\qquad 1\le i \le n.
\end{equation}
Applying $T_\gamma\otimes I$ to (\ref{eq3.5}), and using (\ref{eq3.4}) and the
triangle inequality, gives
\begin{equation}\label{eq3.6}
\|T_\gamma RS_\gamma\otimes I(x_i) - x_i\| < \vp, \qquad 1\le i \le n.
\end{equation}
This shows that $\cl{A}$ has the strong OAP.

\n (vii)~~Let $\cl{J}$ be a norm closed ideal in a $C^*$-algebra $\cl{C}$ with
quotient map $\pi$, and suppose that $\cl{B}$ is exact. To show that
$\cl{A}$ is exact, we need only prove that the kernel of
$I_{\cl{A}}~\otimes~\pi$ is contained in $\cl{A} \otimes \cl{J}$, \cite{W2}.
Consider an element $x \in {\mathrm {ker}}(I_{\cl{A}}\otimes \pi)$, and observe
that $S_{\gamma}\otimes I(x) \in {\mathrm{ker}}(I_{{\cl{B}}}\otimes \pi)$,
which is ${\cl{B}}\otimes {\cl{J}}$ by exactness of ${\cl {B}}$. Apply
$T_{\gamma}\otimes I_{\cl{C}}$ and take a limit over $\gamma \in \Gamma$ to
obtain $x \in {\cl{A}}\otimes {\cl{J}}$. This  shows that ${\cl{A}}$ is
exact.   \end{proof}

We close this section by exhibiting one pair of $C^*$-algebras with the CCFP.

\begin{thm}\label{thm3.2}
Let $G$ be a locally compact group and let $\alpha\colon \ G\to \text{Aut}(\cl{A})$
be a strongly continuous action on a $C^*$-algebra $\cl{A}$. Then $(\cl{A},
\cl{A}\times_{\alpha,r}G)$ has the CCFP.
\end{thm}

\begin{proof}
Consider $a_1, \ldots, a_n\in \cl{A}$ and $\vp>0$. We fix a positive number $\vp'$
to be chosen later. For each $f\in K(G)$, $a\in \cl{A}$, let $f\cdot a$ denote the
element of $K(G,\cl{A})$ whose value at $s\in G$ is $f(s)a$. Then define $S_f\colon
\ \cl{A}\to \cl{A} \times_{\alpha,r} G$ by
\begin{equation}\label{eq3.7}
S_f(a) = (\tilde\pi\times\lambda)(f\cdot a),\qquad a\in \cl{A}.
\end{equation}
Then
\begin{equation}\label{eq3.8}
S_f(a) = \int f(s) \tilde\pi(a)\lambda_s \  ds.
\end{equation}
Each map $a\to \tilde\pi(a)\lambda_s$ is a complete contraction, so $S_f$ may
be viewed as a vector integral of such maps, giving $\|S_f\|_{cb}\le \|f\|_1$.

For each $\xi\in L^2(G)\cap K(G)$, let $\omega_\xi$ be the associated normal
vector functional on $B(L^2(G))$. Then the left slice map $L_{\omega_\xi}$ is
well defined on $B(H) \overline \otimes B(L^2(G))$, (which we identify with
$B(L^2(G,H))$, and $\|L_{\omega_\xi}\|_{cb} = \|\xi\|^2_2$. Let $\widetilde
T_\xi$ be the restriction of $L_{\omega_\xi}$ to $\cl{A}\times_{\alpha,r}G$. Then
$\|\widetilde T_\xi\|_{cb} \le \|\xi\|^2_2$. We first show that the range of
$\widetilde T_\xi$ is contained in $\pi(\cl{A})$. In the following calculation, our
assumptions on continuity of $\xi$ and $f$ will automatically imply that the
integrands are integrable. If $h,k\in H$, then
\begin{align}
\langle\widetilde T_\xi(\tilde\pi\times\lambda)(f\cdot a) h,k\rangle &=
\langle(\tilde\pi \times\lambda)(f\cdot a) h\otimes\xi, k\otimes \xi\rangle
\nonumber\\
&= \iint \langle f(s)\tilde\pi(a)\lambda_s(h\otimes \xi)(t), k\otimes
\xi(t)\rangle dsdt\nonumber\\
&= \iint f(s) \xi(s^{-1}t) \overline{\xi(t)} \langle \pi(\alpha_{t^{-1}}
(a)h,k\rangle dsdt\nonumber\\
\label{eq3.9}
&= \int f*\xi(t) \overline{\xi(t)} \langle\pi(\alpha_{t^{-1}}(a) h,k\rangle dt.
\end{align}
It follows from (\ref{eq3.9}) that
\begin{equation}\label{eq3.10}
\widetilde T_\xi(\tilde\pi\times\lambda)(f\cdot a) = \int f*\xi(t)
\overline{\xi(t)} \pi(\alpha_{t^{-1}}(a))dt
\end{equation}
and this last integral is an element of $\pi(\cl{A})$. If we let $T_\xi$ denote
$\pi^{-1}\widetilde T_\xi$, then we have shown that $\|T_\xi\|_{cb} \le
\|\xi\|^2_2$, and $T_\xi$ maps $\cl{A}\times_{\alpha,r}G$ into $\cl{A}$, since the span
of the elements $\tilde\pi\times\lambda(f\cdot a)$, $f\in K(G)$, $a \in {\cl
A}$, is norm dense in $\cl{A}\times_{\alpha,r}G$.

For the given elements $a_1,\ldots, a_n$ we now wish to choose $f$ and $\xi$
so that
\begin{equation}\label{eq3.11}
\|T_\xi S_f(a_i) - a_i\| < \vp,\qquad 1\le i \le n.
\end{equation}
We restrict attention to $f\in K(G)^+$, $\|f\|_1 = 1$, and $\xi\in L^2(G)^+$,
$\|\xi\|_2 =1$, so that we already have $\|T_\xi\|_{cb}$, $\|S_f\|_{cb}\le 1$.
For each $i$, the map $t\to \alpha_t(a_i)$ is continuous on $G$, so we may
choose a symmetric neighborhood $U$ of $e\in G$ such that
\begin{equation}\label{eq3.12}
\|\alpha_t(a_i)-a_i\| <\vp',\qquad  1\le i \le n,\quad t\in U.
\end{equation}
Now choose a non-negative $\xi\in L^2(G)\cap K(G)$, $\|\xi\|_2 =1$, whose
support is contained in $U$, and let ${}_s\xi$ denote the left translate
$\xi(s^{-1}t)$ of $\xi$. The map $s\to \|\xi-{}_s\xi\|_2$ is continuous on
$G$, by Theorem~20.4 of \cite{HR}, so there is a neighborhood $V$ of $e$,
contained in $U$, within which $\|\xi -{}_s\xi\|_2 < \vp'$. In particular, the
Cauchy-Schwarz inequality shows that \begin{equation}\label{eq3.13}
|\langle{}_s\xi-\xi,\xi\rangle|<\vp',\qquad s\in V.
\end{equation}
Finally we choose $f\in K(G)^+$, $\|f\|_1=1$, and having support in $V$. We
are now ready to show that, with these choices, $\widetilde T_\xi S_f(a_i)$ is
close to $\pi(a_i)$ for $1\le i \le n$.

Let $h,k\in H$ be arbitrary vectors of unit norm. Then
\begin{align}
|\langle(\widetilde T_\xi S_f(a_i) - \pi(a_i)) h,k\rangle|
&= \left|\iint f(s) \xi(s^{-1}t) \xi(t) \langle\pi(\alpha_{t^{-1}}(a_i))
h,k\rangle dsdt - \langle\pi(a_i)h,k\rangle\right|\nonumber\\
&\le \left|\iint f(s) \xi(s^{-1}t) \xi(t) \langle\pi(a_i) h,k\rangle
dsdt - \langle \pi(a_i) h,k\rangle\right|\nonumber\\
\label{eq3.14}
&\quad + \left|\iint f(s) \xi(s^{-1}t) \xi(t) \langle\pi(\alpha_{t^{-1}}
(a_i) - a_i) h,k\rangle dsdt\right|.
\end{align}
We now estimate these integrals separately. Since
\begin{equation}\label{eq3.15}
\iint f(s) \xi(t)^2 dsdt = 1,
\end{equation}
the first may be rewritten as
\begin{align}
\left|\iint \langle\pi(a_i) h,k\rangle f(s) \xi(t) ({}_s\xi(t)-\xi(t))
dsdt\right|
&\le \int_V \|a_i\| f(s)|\langle {}_s\xi-\xi,\xi\rangle| ds\nonumber\\
&\le \int \vp'\|a_i\| f(s) \ ds\nonumber\\
\label{eq3.16}
&= \vp'\|a_i\|.
\end{align}
Here we have used (\ref{eq3.13}) and Fubini's theorem, which is permissible
because the integrand is a continuous function of compact support on $G\times 
G$.

For the second integral, we change the order of integration, and observe that
the $t$-variable can be restricted to $U$ (because of the $\xi(t)$ term). This
yields
\begin{align}
\left|\iint f(s) \xi(s^{-1}t) \xi(t) \langle\pi(\alpha_{t^{-1}} (a_i) -
a_i) h,k\rangle dtds\right|
&\le \iint f(s) \xi(s^{-1}t) \xi(t) \|\alpha_{t^{-1}}(a_i) - a_i\|
dtds\nonumber\\
&\le \iint \vp' f(s) \xi(s^{-1}t) \xi(t) dtds\nonumber\\
&= \int \vp'f(s) [\langle{}_s\xi-\xi,\xi\rangle + \langle
\xi,\xi\rangle]ds\nonumber\\
\label{eq3.17}
&\le (\vp')^2 + \vp',
\end{align}
again using (\ref{eq3.13}). Returning to (\ref{eq3.14}), the two estimates
(\ref{eq3.16}) and (\ref{eq3.17}) lead to
\begin{equation}\label{eq3.18}
|\langle\widetilde T_\xi S_f(a_i) - \pi(a_i) h,k\rangle| \le \vp'\|a_i\| +
\vp' + (\vp')^2.
\end{equation}
The unit vectors were arbitrary in (\ref{eq3.18}), so
\begin{equation}\label{eq3.19}
\|T_\xi S_f(a_i) - a_i\| \le \vp'(1 + \|a_i\| + \vp').
\end{equation}
The proof is completed by choosing $\vp'$ so small  that
\begin{equation}\label{eq3.20}
\vp'(1+\vp' + \max_i \|a_i\|)<\vp.
\end{equation}
Thus $(\cl{A}, \cl{A}\times_{\alpha,r}G)$ has the CCFP.
\end{proof}

\begin{cor}\label{cor3.3}
For any $C^*$-dynamical system $(\cl{A},G,\alpha)$,
\begin{equation}\label{eq3.21}
\Lambda(\cl{A}) \le \Lambda(\cl{A}\times_{\alpha,r}G).
\end{equation}
\end{cor}

\begin{proof}
Apply Theorem~\ref{thm3.1} (3) (i) to the pair $(\cl{A},\cl{A}\times_{\alpha,r}G)$,
which has the CCFP by Theorem~\ref{thm3.2}.
\end{proof}

In \cite{SS2}, the equality $\Lambda(\cl{A}) = \Lambda(\cl{A} \times_{\alpha,r}G)$ was
proved when $G$ was amenable and discrete, or abelian and compact.
Theorem~\ref{thm3.2} already allows us to improve the situation.

\begin{cor}\label{cor3.4}
Let $G$ be an abelian locally compact group, and let $\alpha$ be an action on
a $C^*$-algebra $\cl{A}$. Then
\begin{equation}\label{eq3.22}
\Lambda(\cl{A}) = \Lambda(\cl{A}\times_{\alpha,r}G).
\end{equation}
\end{cor}

\begin{proof}
Since $G$ is abelian, the Takai duality theorem, \cite{T}, states that there
is a dual action $\hat\alpha$ of $\widehat G$ on $\cl{A}\times_{\alpha,r}G$, and
$(\cl{A}\times_{\alpha,r} G) \times_{\hat\alpha,r} \widehat G$ is isomorphic to
$\cl{A}\otimes K(L^2(G))$. Two applications of Corollary~\ref{cor3.3} give
\begin{align}
\Lambda(\cl{A}) &\le \Lambda(\cl{A}\times_{\alpha,r} G) \le \Lambda((\cl{A}\times_{\alpha,r}
G) \times_{\hat\alpha,r} \widehat G)\nonumber\\
&= \Lambda(\cl{A}\otimes K(L^2(G))\nonumber\\
\label{eq3.23}
&= \Lambda(\cl{A}),
\end{align}
the last equality being a special case of Theorem~2.2 of \cite{SS1}.
\end{proof}

\section{Coactions}\label{sec4}

\setcounter{equation}{0}

\indent

In this section we investigate the counterparts of the results of the
previous section for coactions of groups on $C^*$-algebras. One might hope
that the pair $(\cl{A}, \cl{A}\times_\delta G)$ has the CCFP, but this is not possible
in general (we return to this point later). However, if $\cl{A}$ is replaced by a
reduced crossed product $\cl{A}\times_{\alpha,r} G$, the dual coaction $\hat\alpha$
is easy to describe, and a concrete faithful representation of
$(\cl{A}\times_{\alpha,r}G)\times_{\hat\alpha} G$ can be given on the Hilbert space
$H\otimes L^2(G) \otimes L^2(G)$, when $\cl{A}$ is faithfully represented on $H$ by
$\pi$, \cite{IT}. As usual, we identify this Hilbert space with $L^2(G\times G,
H)$ or $L^2(G, L^2(G,H))$ whenever convenient. It is easiest  to work with the
elements $\tilde\pi \times\lambda(f\cdot a)$, $f\in K(G)$, $a\in \cl{A}$, and since
these span a dense subspace of $\cl{A}\times_{\alpha,r}G$, nothing is lost by
making this restriction. The algebra $(\cl{A}\times_{\alpha,r}G)\times_{\hat\alpha}
G$ is generated by a copy of $\cl{A}\times_{\alpha,r}G$ and a copy
 of $C_0(G)$ as multiplication operators on the
second copy of $L^2(G)$. Then, for $a\in \cl{A}$, $h_1,h_2\in H$, $\xi_1$, $\xi_2$,
$\eta_1,\eta_2\in L^2(G) \cap K(G)$, $f,g\in K(G)$, a routine calculation
shows  that
\begin{align}
&\langle \hat\alpha(\tilde\pi\times \lambda(f\cdot a)) I\otimes I\otimes M_g
(h_1\otimes \xi_1 \otimes \eta_1), h_2\otimes \xi_2\otimes
\eta_2\rangle\nonumber\\ \label{eq4.1}
&\quad = \iiint f(s) \langle\pi(\alpha_{r^{-1}}(a)) h_1,h_2\rangle
\xi_1(s^{-1}r)\eta_1(s^{-1}t) g(s^{-1}t)
 \overline{\xi_2(r)} \overline{\eta_2(t)} drdsdt.
\end{align}
Our requirement that all the functions lie in $K(G)$ means that we need not
concern ourselves with the order of integration here, or subsequently.

\begin{thm}\label{thm4.1}
Let $G$ be an amenable locally compact group, let $(\cl{A},G,\alpha)$ be a
$C^*$-dynamical system, and let $\hat\alpha$ be the dual coaction of $G$ on
$\cl{A}\times_{\alpha,r} G$. Then the pair \newline $(\cl{A}\times_{\alpha,r}
G, (\cl{A}\times_{\alpha,r} G) \times_{\hat\alpha} G)$ has the CCFP.
\end{thm}

\begin{proof}
It is clearly sufficient to consider a finite number of elements $\tilde\pi
\times \lambda(f_i\cdot a_i) \in \cl{A}\times_{\alpha,r} G$, $1\le i \le n$, where
$f_i\in K(G)$ and $a_i\in \cl{A}$, $\|f_i\|_1 = 1$, $\|a_i\|=1$. For $g\in K(G)$,
define $S_g\colon \ \cl{A}\times_{\alpha,r}G\to  (\cl{A}\times_{\alpha,r}
G)\times_{\hat\alpha} G$ by
\begin{equation}\label{eq4.2}
S_g(x) =  \hat\alpha(x) I\otimes I\otimes M_g,\qquad x\in \cl{A}\times_{\alpha,r}G.
\end{equation}
Since $\hat\alpha$ is a homomorphism, it is clear that $\|S_g\|_{cb} \le
\|g\|_\infty$. For each $\eta\in L^2(G) \cap K(G)$, we define a map
$T_\eta\colon \ (\cl{A}\times_{\alpha,r}G)\times_{\hat\alpha} G\to
\cl{A}\times_{\alpha,r}G$, $\|T_\eta\|_{cb}\le \|\eta\|^2_2$, by slicing in the
second copy of $L^2(G)$ by the vector functional $\omega_\eta$. {\it A
priori}, $T_\eta$ maps into $B(H\otimes L^2(G))$, but it will become apparent
from the subsequent calculations that the range of $T_\eta$ lies in
$\cl{A}\times_{\alpha,r}G$. From (\ref{eq4.1}),
\begin{align}
&\langle T_\eta(S_g(\tilde\pi \times \lambda(f\cdot a))) h_1 \otimes \xi_1, h_2
\otimes \xi_2\rangle\nonumber\\
\label{eq4.2a}
&\quad = \iiint f(s) \langle \pi(\alpha_{r^{-1}}(a) h_1,h_2\rangle
\xi_1(s^{-1}r) \overline{\xi_2(r)}
\eta(s^{-1}t) g(s^{-1}t) \overline{\eta(t)} drdsdt.
\end{align}
Let $F_{\eta,g} \in K(G)$ be defined by
\begin{equation}\label{eq4.3}
F_{\eta,g}(s) = \int \eta(s^{-1}t) g(s^{-1}t) \overline{\eta(t)} dt,\qquad
s\in G.
\end{equation}
Then we may integrate first with respect to $t$ in (\ref{eq4.2a}) to conclude
that
\begin{equation}\label{eq4.4}
T_\eta(S_g(\tilde\pi \times \lambda(f\cdot a))) = \tilde\pi\times \lambda
((F_{\eta,g}f)\cdot a).
\end{equation}
This establishes that each $T_\eta$ maps into $\cl{A}\times_{\alpha,r}G$.  It also
follows from (\ref{eq4.4}) that
\begin{equation}\label{eq4.5}
T_\eta(S_g(\tilde\pi \times\lambda (f_i\cdot a_i))) - \tilde\pi
\times\lambda(f_i\cdot a_i) = \tilde\pi \times \lambda(((F_{\eta,g}
-1)f_i)\cdot a_i), \end{equation}
so to show that $T_\eta S_g$ is approximately the identity on the elements
$\tilde\pi \times\lambda(f_i\cdot a_i)$, it suffices to choose $\eta$ and $g$
so that the right-hand side of (\ref{eq4.5}) is small in norm. A simple
estimate gives
\begin{equation}\label{eq4.6}
\|\tilde\pi\times \lambda((F_{\eta,g}-1) f_i\cdot a_i)\| \le \|(F_{\eta,g}-1)
f_i\|_1,
\end{equation}
since $\|a_i\| =1$, so given $\vp>0$, it suffices to find $\eta\in L^2(G)
\cap K(G)$, $\|\eta\|_2 \le 1$, $g\in C_0(G)$, $\|g\|_\infty\le 1$, such that
\begin{equation}\label{eq4.7}
|F_{\eta,g}-1|   < \vp
\end{equation}
on the combined supports of the $f_i$'s.

We now use the hypothesis that $G$ is amenable. Let $E_1$, a compact subset of
$G$, be the union of the supports of the $f_i$'s. By Proposition~7.3.8 of
\cite{P}, there is a unit vector $\eta\in L^2(G)$ (which we may take to be in
$L^2(G) \cap K(G))$ such that
\begin{equation}\label{eq4.8}
\eta *\tilde\eta(s) = \int \eta(t) \overline{\eta(s^{-1}t)} dt
\end{equation}
satisfies
\begin{equation}\label{eq4.9}
|\eta*\tilde\eta(s) - 1| < \vp, \qquad s\in E_1.
\end{equation}
Here $\tilde \eta$ is defined by $\tilde\eta(s) = \overline{\eta(s^{-1})}$.
Let $E_2$ be the support of $\eta$, and select $g\in K(G)$, $\|g\|_\infty =
1$, and $g\equiv 1$ on $E_2$. For these choices, it follows from
(\ref{eq4.3}), (\ref{eq4.8}) and (\ref{eq4.9}) that
\begin{equation}\label{eq4.10}
|F_{\eta,g}(s)-1| <\vp, \qquad s\in E_1.
\end{equation}
An immediate consequence of (\ref{eq4.10}) is 
\begin{equation}\label{eq4.11}
\|(F_{\eta,g}-1)f_i\|_1 < \vp, \qquad 1\le i \le n,
\end{equation}
completing the proof.
\end{proof}

We can now state one of our main results.

\begin{cor}\label{cor4.2}
If $G$ is an amenable group and $(\cl{A},G,\alpha)$ is a $C^*$-dynamical system,
then
\begin{equation}\label{eq4.12}
\Lambda(\cl{A}) = \Lambda(\cl{A}\times_{\alpha,r}G).
\end{equation}
\end{cor}

\begin{proof}
We have already shown that $\Lambda(\cl{A}) \le \Lambda(\cl{A}\times_{\alpha,r}
G)$. For the converse, let $\hat\alpha$ be the dual coaction. By
Theorem~\ref{thm4.1}, the pair $(\cl{A}\times_{\alpha,r}G,
(\cl{A}\times_{\alpha,r}G)\times_{\hat\alpha} G)$ has the CCFP, so by
Theorem~\ref{thm3.1} (3) (i) \begin{equation}\label{eq4.13}
\Lambda(\cl{A}\times_{\alpha,r}G) \le \Lambda((\cl{A}\times_{\alpha,r}G)
\times_{\hat\alpha} G).
\end{equation}
But $(\cl{A}\times_{\alpha,r}G)\times_{\hat\alpha} G$ is isomorphic to $\cl{A}\otimes
K(L^2(G))$, \cite{IT}, so (\ref{eq4.13}) becomes
\begin{equation}\label{eq4.14}
\Lambda(\cl{A}\times_{\alpha,r} G) \le \Lambda(\cl{A}\otimes K(L^2(G))) = \Lambda(\cl{A})
\end{equation}
by Theorem 2.2 of \cite{SS1}, proving (\ref{eq4.12}).
\end{proof}

We may now use non-abelian duality to obtain the counterpart of
Corollary~\ref{cor4.2} for coactions. Curiously, it does not seem possible to
prove this without first considering the special case of Theorem~\ref{thm4.1}.
A similar situation arose in the fourth section of \cite{GL}.
\begin{cor}\label{cor4.3}
Let $\delta$ be a non-degenerate coaction of an amenable group $G$ on a
$C^*$-algebra $\cl{A}$. Then
\begin{equation}\label{eq4.15}
\Lambda(\cl{A}\times_\delta G) = \Lambda(\cl{A}).
\end{equation}
\end{cor}

\begin{proof}
By \cite{Ka}, there is a dual action $\hat\delta$ of $G$ on $\cl{A}\times_\delta G$
such that $(\cl{A}\times_\delta G)\times_{\hat\delta,r}G$ is isomorphic to
$\cl{A}\otimes K(L^2(G))$. Using Corollary~\ref{cor4.2}, with $\cl{A}\times_\delta G$ in
place of $\cl{A}$, we see that
\begin{equation}\label{eq4.16}
\Lambda(\cl{A}\times_\delta G) = \Lambda((\cl{A}\times_\delta G)
\times_{\hat\delta,r} G) = \Lambda(\cl{A} \otimes K(L^2(G))) = \Lambda(\cl{A}),
\end{equation}
proving (\ref{eq4.15}).\end{proof}

\begin{rem}\label{rem4.4}
If $G$ is non-amenable then (\ref{eq4.15}) may fail, but the inequality
\begin{equation}\label{eq4.17}
\Lambda(\cl{A}\times_\delta G) \le \Lambda(\cl{A})
\end{equation}
is an immediate consequence of applying our duality methods and
Corollary~\ref{cor3.3}.$\hfill\square$
\end{rem}

\begin{rem}\label{rem4.5}
Corollary~\ref{cor4.2} also holds true for the twisted crossed products of
\cite{PR}. When $G$ is amenable, full and reduced crossed products coincide,
so from \cite{PR}, a twisted crossed product of $\cl{A}$ by $G$ is stably
isomorphic to a crossed product of $\cl{A}$ by $G$. Our assertion then follows
from Corollary~\ref{cor4.2}, since tensoring by the algebra of compact
operators does not change $\Lambda(\cdot)$.$\hfill\square$
\end{rem}

We conclude this section by investigating whether actions and coactions
preserve the properties stated in Theorem~\ref{thm3.1}.

\begin{thm}\label{thm4.6}
Let $G$ be an amenable group and let $\alpha,\delta$ be respectively an action
and a coaction of $G$ on a $C^*$-algebra $\cl{A}$. For any one of the properties
(i)--(vii) of Theorem~\ref{thm3.1}, all three $C^*$-algebras $\cl{A},
\cl{A}\times_{\alpha,r}G$, and $\cl{A}\times_\delta G$ have this property or all three
do not.
\end{thm}

\begin{proof}
Since the pair $(\cl{A},\cl{A}\otimes K(H))$ has the CCFP for any Hilbert space $H$
(Theorem~\ref{thm3.1} (1)), these properties all transfer from $\cl{A}\otimes K(H)$
to $\cl{A}$. The verification that these properties all transfer from $\cl{A}$ to
$\cl{A}\otimes K(H)$ is essentially routine, based on the nuclearity of $K(H)$. To
give the flavor, we will prove this for (iii), and leave the others to the
reader.

Suppose that $\cl{A}$ has the general slice map property and fix a $C^*$-algebra
$\cl{B}$ with a closed subspace $\cl{E}$. Let $x\in (\cl{A}\otimes K(H))\otimes \cl{B}$ be an
element whose right slices lie in $\cl{E}$. If $\phi\in \cl{A}^*$ and $\psi\in K(H)^*$,
then slicing by $\phi\otimes \psi$ is the same as slicing by $\phi$ and then
by $\psi$. Thus the element $y\in K(H)\otimes \cl{B}$ obtained from $x$ by slicing
by $\phi$ has the property that all right slices are in $\cl{E}$. Since $K(H)$ has
the general slice map property, by nuclearity, it follows that $y\in K(H)
\otimes \cl{E}$. Since $\phi\in \cl{A}^*$ was arbitrary, the general slice map property
for $\cl{A}$ implies that $x\in \cl{A}\otimes K(H)\otimes \cl{E}$ as required.

We have already shown in Theorems~\ref{thm3.2} and \ref{thm4.1} that the pairs
$(\cl{A},\cl{A}\times_{\alpha,r}G)$ and \newline
$(\cl{A}\times_{\alpha,r}G,
\cl{A}\otimes K(L^2(G))$ have the CCFP, so we conclude from
Theorem~\ref{thm3.1} and the preceding remarks that these properties transfer
between $\cl{A}$ and $\cl{A}\times_{\alpha,r} G$ in both directions. Applying
this to the pair $\cl{A}\times_\delta G$ and $(\cl{A}\times_\delta G)
\times_{\hat \delta,r}G \approx \cl{A}\otimes K(L^2(G))$, we see that they
also transfer in both directions between $\cl{A}$ and $\cl{A}\times_\delta G$.
\end{proof}

\begin{rem}\label{rem4.7}
We note that (ii) in Theorem~\ref{thm4.6} gives a new method of showing the
well known result (see \cite{G,LP,R}) that nuclearity is preserved by actions
and coactions of amenable groups. The same is also true for
exactness, by (vii), where the original proofs of the action and coaction cases
are due respectively to Kirchberg and to Ng, \cite{Kir2,Ng}.$\hfill\square$
\end{rem} \newpage \section{Concluding remarks}\label{sec5}

\indent

If $\alpha$ is the trivial action of $G$ on $\cl{A}$, then the crossed product
$\cl{A}\times_{\alpha,r}G$ is isomorphic to $\cl{A}\otimes C^*_r(G)$. It then follows
from Theorem~2.2 of \cite{SS1} that $\Lambda(\cl{A}\times_{\alpha,r} G) =
\Lambda(\cl{A}) \Lambda(C^*_r(G))$ in this case. One might hope that this formula
holds in general or, failing this, that $\Lambda(\cl{A}\times_{\alpha,r}G)$ is
bounded above by an expression involving only $\Lambda(\cl{A})$ and
$\Lambda(C^*_r(G))$. The following example, which comes from \cite{H,HK} shows
that this is impossible.

Let a group $G$ act on another group $H$ by a homomorphism $\rho\colon \ G\to
\text{Aut}(H)$ where both are discrete. There is a naturally induced action
$\alpha$ of $G$ on $C^*_r(H)$, and the crossed product $C^*_r(H)
\times_{\alpha,r}G$ is isomorphic to the reduced $C^*$-algebra
$C^*_r(H\times_\rho G)$ of the semi-direct product $H\times_\rho G$, \cite{B}.
It was established in \cite{H} that, for the natural action of
$SL(2,\mathbb{Z})$ on $\mathbb{Z}^2$, $\Lambda (C^*_r(\mathbb{Z}^2\times_\rho
SL(2,\mathbb{Z}))) = \infty$, while $\Lambda(C^*_r(\mathbb{Z}^2))$ and
$\Lambda(C^*_r(SL(2,\mathbb{Z}))$ are both finite. Thus we cannot expect an
upper estimate for $\Lambda(\cl{A}\times_{\alpha,r}G)$ in terms of $\Lambda(\cl{A})$ and
$\Lambda(C^*_r(G))$.

Any group $G$ acts trivially on ${\mathbb C}$, and the resulting reduced
crossed product is $C^*_r(G)$. By duality, $C^*_r(G)\times_{\delta}G$
(where $\delta$ is the dual coaction) is isomorphic to $K(H)$. Thus the
pair $(C^*_r(G),C^*_r(G)\times_{\delta}G)$ can only have the CCFP when
$C^*_r(G)$ is nuclear, confirming our statement at the beginning of the
previous section that pairs $({\cl A},{\cl A}\times_{\delta}G)$ will not always
have this property.

In \cite{HK}, Haagerup and Kraus introduced the approximation property (AP)
for groups (we refer to this paper for the definition which will not be needed
here). For a discrete group $\Gamma$, Theorem~2.1 of \cite{HK} shows that the
AP for $\Gamma$ is equivalent to the strong OAP for $C^*_r(\Gamma)$. They left
open the question of whether the AP passes to quotients by normal subgroups,
but expressed the view that this was unlikely since every countable discrete
group is a quotient of $\mathbb{F}_\infty$, which {\it does\/} have the AP,
\cite{DH}. Since the strong OAP implies exactness (remarks preceding
Theorem~2.2 of \cite{HK}), the work of Ozawa, \cite{O}, showing the existence
of discrete groups $\Gamma$ for which $C^*_r(\Gamma)$ is not exact, also
provides examples where the AP fails. The above discussion then shows that
there are normal subgroups $N$ of $\mathbb{F}_\infty$ such that
$C^*_r(\mathbb{F}_\infty)$ and $C^*_r(N)$ have the strong OAP while the
$C^*$-algebras $C^*_r(\mathbb{F}_\infty/N)$ do not.

Although we have not needed to do so, we could have weakened the definition
of the CCFP by requiring that the nets of maps be uniformly bounded in the
completely bounded norm. The {\it completely bounded factorization
property} (CBFP) would be an appropriate name for this potentially useful
concept. Most of the statements and proofs of Theorems~\ref{thm3.1} and
\ref{thm4.6}  are valid with little change, although for {\it nuclearity} a
result from \cite{Pi} (Remark preceding Theorem 2.10) is
required: ${\cl A}$ is nuclear if and only if $({\cl A},K(H))$ has the
CBFP.

\end{document}